\documentclass[12pt]{amsart}
\usepackage{amssymb}
\usepackage{latexsym}

\title{Line partitions of internal points to a conic in $PG(2,q)$}
\author{Massimo Giulietti}
\date{}

\newtheorem{proposition}{Proposition}[section]
\newtheorem{theorem}[proposition]{Theorem}

\newtheorem{lemma}[proposition]{Lemma}

\newcommand{\ged}{\vspace{.3cm}\hfill$\blacksquare$}
\def\cA{{\mathcal A}}
\def\ki{\chi}

\def\cI{{\mathcal I}}
\def\cF{{\mathcal F}}
\def\cR{{\mathcal R}}

\def\cC{\mathrm C}

\def\cL{\mathcal L}

\def\cF{\mathcal F}

\def\fq{{\mathbb F}_q}

\begin{document}

\begin{abstract} All sets of lines providing a partition of the set of internal
points to a conic $\cC$ in $PG(2,q)$, $q$ odd, are determined.
There exist only three such linesets up to projectivities, namely
the set of all nontangent lines to $\cC$ through an external point
to $\cC$, the set of all nontangent lines to $\cC$ through a point
in $\cC$, and, for square $q$, the set of all nontangent lines to
$\cC$ belonging to a Baer subplane $PG(2,\sqrt{q})$ with
$\sqrt{q}+1$ common points with $\cC$. This classification theorem
is the analogous of a classical result by Segre and Korchm\'aros
\cite{KS} characterizing the pencil of lines through an internal
point to $\cC$ as the unique set of lines, up to projectivities,
which provides a partition of the set of all noninternal points to
$\cC$. However, the proof is not analogous, since it does not rely
on the famous Lemma of Tangents of Segre which was the main
ingredient in \cite{KS}. The main tools in the present paper are
certain partitions in conics of the set of all internal points to
$\cC$, together with some recent combinatorial characterizations
of blocking sets of non-secant lines, see \cite{AK2}, and of
blocking sets of external lines, see \cite{AK1}.
\end{abstract}

\maketitle

\section{Introduction} In 1977 Segre and Korchm\'aros gave the
following combinatorial characterization of external lines to an
irreducible conic in $PG(2,q)$, see \cite{KS}, \cite{H} Theorem
13.40, and \cite{KK}.
\begin{theorem}
\label{segrekorchmaros} If every secant and tangent of an
irreducible conic meets a pointset $\cL$ in exactly one point, then
$\cL$ is linear, that is, it consists of all points of an external
line to the conic.
\end{theorem}
For even $q$, this was proven by Bruen and Thas \cite{BT},
independently.

It is natural to ask for a similar characterization of a minimal
pointset $\cL$ meeting every {\em external line} to an irreducible
conic $\cC$ in exactly one point. In this case, we have two linear
examples: a chord minus the common points with $\cC$, and a
tangent minus the tangency point (and, for $q$ even, minus the
nucleus of $\cC$, as well).

For $q$ even, it is shown in \cite{MG} that there is exactly one
more possibility for $\cL$, namely, for any even square $q$, the
set consisting of the points of a Baer subplane $\pi$ sharing
$\sqrt{q}+1$ with $\cC$, minus $\pi\cap \cC$ and the nucleus of
$\cC$.

The aim of the present paper is to prove an analogous result for
$q$ odd.

Henceforth, $q$ is always assumed to be odd, that is, $q=p^h$ with
$p>2$ prime. Then the orthogonal polarity associated to $\cC$
turns $\cL$ into a {\em line partition} of the set of all internal
points to $\cC$. In terms of a line partition, Theorem
\ref{segrekorchmaros} states that if $\cL$ is a line partition of
the set of all noninternal points to $\cC$, then $\cL$ is a pencil
of lines through an internal point to $\cC$.

Our main result is the following theorem.
\begin{theorem}\label{main}
Let $\cL$ be a line partition of the set of internal points to a
conic $\cC$ in $PG(2,q)$, $q$ odd. Then either
\begin{itemize}
\item $\#\cL=q-1$, and $\cL$ consists of the $q-1$ lines through
an external point of $\cC$ which are not tangent to $\cC$, or

\item $\#\cL=q$, and $\cL$ consists of the $q$ lines through a
point of $\cC$ distinct from the tangent to $\cC$, or

\item $\#\cL=q$ for a square $q$, and $\cL$ consists of all non
tangent lines belonging to a Baer subplane $PG(2,\sqrt q)$ with
$\sqrt{q}+1$ common points with $\cC$.
\end{itemize}
\end{theorem}
%
%
%

\section{Internal points to a conic}
In this section a certain partition in conics of the internal
points to a conic $\cC$ in $PG(2,q)$, $q$ odd, is investigated.

Assume without loss of generality that $\cC$ has affine equation
$Y=X^2$, and denote by $Y_\infty$ the infinite point of $\cC$.
Consider the pencil of conics $\cF$ consisting of the conics
$\cC_s:Y=X^2-s$, with $s$ ranging over $\fq$.

First, an elementary property of $\cF$ which will be useful in the
sequel is pointed out.
\begin{lemma}\label{oss2}
Any line of $PG(2,q)$ not passing through $Y_\infty$ is tangent to
exactly one conic of $\cF$.
\end{lemma}
{\em Proof.} It is enough to note that the line of equation
$Y=\alpha X+ \beta$ is tangent to $\cC_s$ if and only if
$s=-\frac{\alpha^2+4\beta}{4}$.\ged

Recall that in the finite field $\fq$ half the non-zero elements
are quadratic residues or squares, and half are quadratic
non-residues or non-squares. The quadratic character of $\fq$ is
the function $\ki$ given by
$$
\ki(x)=\left\{\begin{array}{rl}0& \text{if
}x=0\,,\\
1& \text{if }x \text{ is a quadratic residue,}\\
-1& \text{if }x \text{ is a quadratic non-residue.}
\end{array}\right.
$$

\begin{lemma}\label{l1}
Let $\cC_{s}$ and $\cC_{s'}$ be two distinct conics in $\cF$. Then
the affine points of $\cC_{s'}$ are all either  external or
internal to $\cC_{s}$, according to whether $\ki(s'-s)=1$ or
$\ki(s'-s)=-1$.
\end{lemma}
{\em Proof.} Let $P=(a,a^2-s')$ be an affine point of $\cC_{s'}$.
The polar line $l_P$ of $P$ with respect to $\cC_{s}$ has equation
$Y=2aX-a^2+s'-2s$. Then it is straightforward to check that $l_P$
does not meet $\cC_s$ if and only if $s'-s$ is a non-square  in
$\fq$. \ged

As a matter of terminology,  we will say that a conic $\cC_s$ is
{\em internal} ({\em external}) to $\cC_{s'}$ if all the affine
points of $\cC_s$ are internal (external) to $\cC_{s'}$. Let
$\cI=\{\cC_s\mid \ki(s)=-1\}$. Clearly, the set of internal points
to $\cC$ consists of the affine points of the conics in $\cI$.

Throughout the rest of this section we assume that $q\equiv 3
\pmod 4$. Note that this is equivalent to $\ki(-1)=-1$, see
\cite{H}. Then Lemma \ref{l1} yields that $\cC_s$ is internal to
$\cC_{s'}$ if and only if $\cC_{s'}$ is external to $\cC_{s}$.

\begin{lemma}\label{nextern} Let $\cC_{s}$ be a conic in $\cI$.
If $q\equiv 3 \pmod 4$, then there are exactly $\frac{q-3}{4}$
conics in $\cI$ that are internal to $\cC_s$.
\end{lemma}
{\em Proof.} The hypothesis $q\equiv 3 \pmod 4$ yields that for
any $s\in \fq$, $s\neq 0$, there are exactly $\frac{q-3}{4}$
ordered pairs $(u,v)\in \fq \times \fq$ with $s=u-v$ and
$\chi(u)=\chi(v)=1$ (see e.g. \cite[Lemma 1.7]{WALL}). Via the
correspondence $s'=-v$, the number of such pairs equals the number
of $s'\in \fq$ satisfying $\ki(s')=\ki(s'-s)=-1$. Then the
assertion follows from Lemma \ref{l1}. \ged

Denote $\cI_s$ the set of conics of $\cI$ which are internal to
$\cC_s$. The following lemma will be crucial in the proof of
Theorem \ref{main}.
\begin{lemma}\label{matrice}
Let $q\equiv 3 \pmod 4$. Then the any integer function $\varphi$
on $\cI$ such that
\begin{equation}\label{chiave0}
\sum_{\cC_{s'}\in\cI_s}\varphi(\cC_{s'})=
\sum_{\cC_{s'}\in\cI\setminus\cI_s,\cC_{s'}\neq\cC_s}
\varphi(\cC_{s'}),\quad \text{ for any }\cC_s\in \cI
\end{equation}
is constant.
\end{lemma}
{\em Proof.} Let $\{s_1, s_2,\ldots,s_{\frac{q-1}{2}}\}$ be the
set of non-squares in $\fq$, and let $A=(a_{ij})$ be the
$\frac{q-1}{2}\times \frac{q-1}{2}$ matrix given by
$$
a_{ij}=\ki(s_i-s_j)\,.
$$
Then by Lemma \ref{l1}, condition (\ref{chiave0}) is equivalent to
$$
\sum_{\ki(s_i-s_j)=-1}\varphi(\cC_{s_i})= \sum_{\ki(s_i-s_j)=1}
\varphi(\cC_{s_i}),\quad \text{ for any
}j=1,\ldots,\frac{q-1}{2}\,,
$$
that is, the vector
$(\varphi(\cC_{s_1}),\ldots,\varphi(\cC_{s_i}),\ldots,\varphi(\cC_{s_{\frac{q-1}{2}}}))$
belongs to the null space of $A$. Clearly if $\varphi$ is constant
such a condition is fulfilled by Lemma \ref{nextern}.

Then to prove the assertion, it is enough to show that the real
rank of $A$ is at least $\frac{q-1}{2}-1$. As usual, denote
$A_{1,1}$  the matrix obtained from $A$ by dismissing the first
row and the first column column. Note that as the entries of
$A_{1,1}$ are integers, $Det(A_{1,1})\pmod 2$ coincides with
$Det({\tilde A_{1,1}})$, where ${\tilde A_{1,1}}$ is the matrix
over the finite field with $2$ elements obtained from $A_{1,1}$ by
substituting each entry $m_{ij}$ with $m_{ij}\pmod 2$. By
definition of $A$, the entries of ${\tilde A_{1,1}}$ are equal to
$1$, except those in the diagonal which are equal to zero. As
$\frac{q-1}{2}-1$ is even, it is straightforward to check that
${\tilde A_{1,1}}^2$ is the identity matrix, whence $Det({\tilde
A_{1,1}})=Det(A_{1,1})\pmod 2$ is different from $0$. \ged

\section{Proof of Theorem \ref{main}}
Throughout, $\cC$ is an irreducible conic in $PG(2,q)$, $q$ odd,
and $\cL$ is a line partition of the set of internal points to
$\cC$. First, the possible sizes of $\cL$ are determined.

\begin{lemma}\label{easy}
The size of $\cL$ is either $q-1$ or $q$. In the latter case,
$\cL$ consists of $q$ secant lines to $\cC$.
\end{lemma}
{\em Proof.} The number of internal points to a conic is
$q(q-1)/2$, see \cite{H}. Also, a secant line of $\cC$ contains
$(q-1)/2$ internal points of $\cC$, whereas the number of internal
points on an external line is $(q+1)/2$. No internal point belongs
to a tangent to $\cC$. Let $\cL$ consist of $h$ secants together
with $k$ external lines to $\cC$. As $\cL$ is a line partition of
the internal points to $\cC$,
$$
\frac{q(q-1)}{2}=h\frac{q-1}{2}+k\frac{q+1}{2}\,,
$$
that is
$$
q=h+k+\frac{2k}{q-1}\,.
$$
As $\frac{2k}{q-1}$ is an integer, either $k=0$ and $h=q$, or
$k=(q-1)/2=h$. \ged

The classification problem for $\#\cL=q-1$ is solved  via the
characterization of blocking sets of minimal size of the external
lines to a conic, as given in \cite{AK1}. The dual of Theorem 1.1
in \cite{AK1} reads as follows.
\begin{proposition}\label{agkor1}
Let $\cR$ be a lineset  of size $q-1$ such that any internal point
to $\cC$ belongs to some line of $\cR$. If either $q=3$ or $q> 9$
, then $\cR$ consists of the $q-1$ lines through an external point
of $\cC$ which are not tangent to $\cC$. For $q=5,7$ there exists
just one more example, up to projectivities, for which some of the
lines in $\cR$ are external to $\cC$.
\end{proposition}

From now on, assume that $\#\cL=q$. Note that Lemma \ref{easy}
yields that every line of $\cL$ is a secant line of $\cC$. We
first deal with the case $q\equiv 3 \pmod 4$.
\begin{lemma}\label{lem1} Let $\#\cL=q$. If $q\equiv 3 \pmod 4$, then
the number of lines of $\cL$ through any point $P$ of $\cC$ is
$1$, $\frac{q+1}{2}$ or $q$.
\end{lemma}
{\em Proof.} We keep the notation of Section 2. Assume without
loss of generality that $\cC$ has equation $X^2-Y=0$, and that
$P=Y_\infty$. Let $\cL_P$ be the set of lines of $\cL$ passing
through $P$, and set $m=\#\cL_P$. Also, for any $l\in \cL\setminus
\cL_P$, denote $\cC^{(l)}$ the conic of $\cF$ which is tangent to
$l$ according to Lemma \ref{oss2}.

As any secant $l$ of $\cC$ not passing through $P$ contains an odd
number of internal points to $\cC$, the conic $\cC^{(l)}$ belongs
to $\cI$. We claim that for any $\cC_s\in \cI$ and for any $l\in
\cL\setminus \cL_P$, $l$ not tangent to $\cC_s$,
\begin{equation}\label{frase}
\cC_s \text{ is external to } \cC^{(l)}  \text{ if and only if } l
\text{ is a secant of } \cC_s\,.
\end{equation}
Clearly, if $l$ is a secant of $\cC_s$, then both the points of
$\cC_s\cap l$  are external to $\cC^{(l)}$. Therefore $\cC_s$ is
external to $\cC^{(l)}$. To prove the only if part of
(\ref{frase}), note that for any $l\in \cL\setminus \cL_P$ the set
of $\frac{q-1}{2}$ points of $l$ which are internal to $\cC$
consists of one point lying on $\cC^{(l)}$ together with
$\frac{q-3}{4}$ point pairs, each of which contained in a conic of
$\cI$. Taking into account Lemma \ref{nextern}, this means that
$l$ is a secant of all the conics of $\cI$ that are external to
$\cC^{(l)}$.

Now, for any $\cC_s\in \cI$ let $\varphi(\cC_s)$ be the number of
lines of $\cL$ which are tangent to $\cC_s$.  Then,
\begin{equation}\label{chiave}
\sum_{\cC_{s'}\in\cI_s}\varphi(\cC_{s'})=
\sum_{\cC_{s'}\in\cI\setminus\cI_s,\cC_{s'}\neq\cC_s}
\varphi(\cC_{s'}),\quad \text{ for any }\cC_s\in \cI\,.
\end{equation}
In fact, (\ref{frase}) yields that
$\sum_{\cC_{s'}\in\cI_s}\varphi(\cC_{s'})$ equals the number of
lines in $\cL\setminus \cL_P$ which are secants to $\cC_s$, that
is $\frac{q-m-\varphi(\cC_s)}{2}$. As the total number of lines in
$\cL$ which are tangent to a conic of $\cI$ distinct from $\cC_s$
is $q-m-\varphi(\cC_s)$, Equation (\ref{chiave}) follows.  Then by
Lemma \ref{matrice}, $\varphi(\cC_s)$ is an integer which is
independent of $\cC_s$. Denote $t$ such an integer. By Lemma
\ref{oss2},
\begin{equation}\label{somma}
\sum_{\cC_s\in \cI}\varphi(\cC_s)=t\frac{q-1}{2}=q-m\,,
\end{equation}
which implies that either (a) $t=2$, $m=1$, (b) $t=0$, $m=q$, or
(c) $t=1$, $m=\frac{q+1}{2}$. \ged

\begin{lemma}\label{lem2}
Let $\#\cL=q$. If $q\equiv 3 \pmod 4$, then no point of $\cC$
belongs to exactly $\frac{q+1}{2}$ lines of $\cL$.
\end{lemma}
{\em Proof.} We keep the notation of the proof of Lemma
\ref{lem1}. Also, for $Q\in \cC$ let $m_Q$ be the number of lines
of $\cL$ passing through $Q$.

Assume that $m_P=\frac{q+1}{2}$, with $P=Y_\infty$. As $\sum_{Q\in
\cC}m_Q=2q$, Lemma \ref{lem1} yields that there exists another
point ${\bar P}\in \cC$ belonging to exactly $\frac{q+1}{2}$ lines
of $\cL$, and that $m_Q=1$ for any  point $Q\in\cC$, $Q\notin
\{P,{\bar P}\}$. As the projective group of $\cC$ is sharply
$3$-transitive on the points of $\cC$ (see e.g. \cite{H}), we may
assume that ${\bar P}$ coincides with $(0,0)$.

Let $\cA$ be the subset of $\fq\setminus \{0\}$ consisting of the
$\frac{q-1}{2}$ non-zero elements $u$ for which the line $Y=uX$
belongs to $\cL$. Then the lines in $\cL_P$ are those of equation
$X=v$, with $v$ ranging over $\fq \setminus \cA$. Actually,
$\fq\setminus \cA$ coincides with $\{-u\mid u \in\cA\}\cup \{0\}$.
In fact, $u\in \cA$ yields $-u\notin \cA$, otherwise the two lines
of equation $Y= u X$ and $Y=-uX$ would be both lines of $\cL$
tangent to the same conic $\cC_{-{u^2}/{4}}$. By the proof of
Lemma \ref{lem1} this is impossible, as $m_P=\frac{q+1}{2}$ yields
that each conic in $\cI$ has exactly one tangent in $\cL\setminus
\cL_P$.

Then, for any $u_1,u_2\in \cA$, $u_1\neq u_2$, the lines $Y=u_1X$
and $X=-u_2$, as well as the lines $Y=u_2X$ and $X=-u_1$, meet in
an external point to $\cC$, that is
$$
\ki(u_1^2+u_1u_2)=\ki(u_2^2+u_2u_1)=1 \,.
$$
Equivalently, for any $u_1,u_2\in \cA$, $u_1\neq u_2$,
$$
\ki(u_1)\ki(u_1+u_2)=\ki(u_2)\ki(u_1+u_2)=1 \,,
$$
whence all the elements in $\cA$ and all the sums of two distinct
elements in $\cA$ have the same quadratic character. But this is
actually impossible, as $q\equiv 3 \pmod 4$ yields that for any
$u_1\in \fq\setminus \{0\}$, $\epsilon \in \{-1,1\}$, the number
of $u_2\in \fq$ such that $\ki(u_2)=\ki(u_1+u_2)=\epsilon$ is
$\frac{q-3}{4}$ (see e.g. \cite[Lemma 1.7]{WALL}). \ged

\begin{proposition}\label{p2} Let $\#\cL=q$. If $q\equiv 3 \pmod 4$, then
 $\cL$ consists of the $q$ lines
through a point of $\cC$ distinct from the tangent to $\cC$.
\end{proposition}
{\em Proof.}  By Lemmas \ref{lem1} and \ref{lem2} the number $m_P$
of lines of $\cL$ through a given point $P\in \cC$ is either $1$
or $q$. As $\#\cL=q>\frac{q+1}{2}$ it  is impossible that $m_P=1$
for every $P\in \cC$. Then there exists a point $P_0$ with
$m_{P_0}=q$, which proves the assertion. \ged

Assume now that $q\equiv 1 \pmod 4$. We first prove that any line
partition  of size $q$ of the internal points of $\cC$ actually
covers all the points of $\cC$ as well.
\begin{lemma}\label{key}
Let $\#\cL=q$. If $q\equiv 1 \pmod 4$, then any point of $\cC$
belongs to some line of $\cL$.
\end{lemma}
{\em Proof.} We keep the notation of Section 2. Assume that a
point $P\in \cC$ does not belong to any line of $\cL$. Without
loss of generality, let $P=Y_\infty$. Then the $q$ affine points
of any conic $\cC_s\in \cI$ partition into sets $l\cap \cC_s$,
with $l$ ranging over $\cL$. As $q$ is odd, there exists a line
$l_s\in \cL$ which is tangent to $\cC_s$. Any line of $\cL$ has an
even number of internal points to $\cC$, as $(q-1)/2$ is even.
Then some line of $\cL$ must be  tangent to more than one conic of
$\cF$, which is a contradiction to Lemma \ref{oss2}.\ged

To complete our investigation for $q\equiv 1\pmod 4$, the
combinatorial characterization of blocking sets of non-secant
lines to $\cC$, as given in \cite{AK2}, is needed. The dual of
Theorem in \cite{AK2} reads as follows.
\begin{lemma}\label{agkor2}
Let $\cR$ be a lineset  of size $q$ such that any non-external
point to $\cC$ belongs to some line of $\cR$. Then one of the
following occurs.
\begin{itemize}
\item[{\rm (a)}] $\cR$ consists of $q$ lines through a point of
$\cC$ distinct from the tangent to $\cC$,

\item[{\rm (b)}] $\cR$ consists of the lines of a subgeometry
$PG(2,\sqrt q)$ which are not tangent to $\cC$.

\item[{\rm (c)}] $\cR$ consists of the $q-1$ lines through an
external point $P$ to $\cC$ which are not tangent to $\cC$,
together with the polar line of $P$ with respect to $\cC$.

\end{itemize}
\end{lemma}

\begin{proposition}\label{p3}
Let $\#\cL=q$. If $q\equiv 1 \pmod 4$, then $\cL$ consists either
of the $q$ lines through a point of $\cC$ distinct from the
tangent to $\cC$, or of the lines of a subgeometry $PG(2,\sqrt q)$
which are not tangent to $\cC$.
\end{proposition}
{\em Proof.} Lemma \ref{key} yields that $\cL$ satisfies the
hypothesis of Lemma \ref{agkor2}. Actually, (c) of Lemma
\ref{agkor2} cannot occur as in this case not every line of $\cR$
is a secant line to $\cC$. Hence the assertion is
 proved.
\ged

Theorem \ref{main} now follows from Propositions \ref{agkor1},
\ref{p2}, \ref{p3}.

\vspace{1cm}\noindent {\em Author's address}:\\

\vspace{0.5 cm} \noindent Massimo GIULIETTI\\  Dipartimento di
Matematica\\  Universit\`a degli studi di Perugia\\Via  Vanvitelli
1\\ 06123 Perugia (Italy).\\
 E--mail: {\tt
 giuliet@dipmat.unipg.it


\begin{thebibliography}{99}
\bibitem{AK1} A. Aguglia, G. Korchm\'aros, {\em Blocking sets of external lines to a conic in PG(2,q),
 $q$ odd}, Combinatorica (to appear).

\bibitem{AK2} A. Aguglia, G. Korchm\'aros, {\em Blocking sets of nonsecant lines to a conic in PG(2,q),
  $q$ odd}, Journal of  Combinatorial Designs (to appear).

\bibitem{AKS} A. Aguglia, G. Korchm\'aros, A. Siciliano, {\em Minimal covering of all chords
of a conic in $PG(2,q)$, $q$ even}, Bulletin of the Belgian
Mathematical Society, (to appear).

\bibitem{BFK} E. Boros, Z. F\H{u}redi, J. Kahn,  {\em Maximal Intersecting Families
and Affine Regular Polygons in $PG(2,q)$}, Journal of
Combinatorial Theory, Series A {\bf 52}, 1--9 (1989).


\bibitem{BT} A. Bruen and J.A. Thas, Flocks, chains and configurations in finite geometries.
{\em Atti Accad. Naz. Lincei Rend. Cl. Sci. Fis. Mat. Natur.} (8)
{\bf{59}} (1975), 744-748 (1976).


\bibitem{H} J.W.P. Hirschfeld, {\em Projective Geometries over Finite
Fields}, Clarendon Press, Oxford (1998).

\bibitem{MG} M. Giulietti, {\em Blocking sets of external lines to a conic in $PG(2,q)$, $q$ even}, submitted.

\bibitem{KK} G. Korchm\'aros {\em Segre's type theorems}, invited lecture at
the International Conference ``Trends in Geometry, in Memory of
Beniamino Segre'', 7-8 June 2004 Rome, to appear in a special
issue of Rendiconti di Matematica e delle sue applicazioni.

\bibitem{KS} B. Segre, G. Korchm\'aros,   {\em Una propriet\'a degli insiemi di punti di un piano di Galois
caratterizzante quelle formati dai punti delle singole rette
esterne ad una conica}, Atti Accad. Naz. Lincei Rend. Cl. Sci.
Fis. Mat. Natur. (8) {\bf 62}, 613--619 (1977).

\bibitem{WALL} W.D. Wallis, A. Penfold Street, J. Seberry Wallis,
{Combinatorics: Room Squares, Sum-Free Sets, Hadamard Matrices},
volume 292 of {\em Lecture Notes in Math.}, Springer-Verlag,
Berlin-Heidelberg-New York, 1972.
\end{thebibliography}
\end{document}